\documentclass[12pt,english]{article}
\usepackage[T1]{fontenc}
\usepackage[utf8]{inputenc}
\usepackage{amsmath}
\usepackage{amsthm}
\usepackage{setspace}
\makeatletter
\theoremstyle{plain}
\newtheorem{assumption}{\protect\assumptionname}
\theoremstyle{plain}
\newtheorem{thm}{\protect\theoremname}
\theoremstyle{remark}
\newtheorem{rem}{\protect\remarkname}
\theoremstyle{plain}
\newtheorem{lem}{\protect\lemmaname}

\usepackage[letterpaper, left=1in, top=1in, right=1in, bottom=1in,nohead,includefoot, verbose, ignoremp]{geometry}
\usepackage{charter} 
\usepackage{enumerate} 
\usepackage{xy}\xyoption{all} \xyoption{poly} \xyoption{knot}  
\usepackage{latexsym,amssymb,amsmath,amsfonts,graphicx,color,fancyvrb,amsthm,enumerate,natbib,movie15}
\usepackage[pdftex,pagebackref=true]{hyperref}
\usepackage[svgnames,dvipsnames,x11names]{xcolor}
\usepackage{graphicx}
\usepackage{graphicx}
\usepackage{yhmath}
\usepackage{mathdots}
\usepackage{MnSymbol}
\hypersetup{
colorlinks,%
linkcolor=RoyalBlue2,  
urlcolor=Sienna4,   
citecolor=RoyalBlue2  
}

\def\beginmat{ \left( \begin{array} }
\def\endmat{ \end{array} \right) }

\setcitestyle{citesep={,}}

\makeatother

\usepackage{babel}
\providecommand{\assumptionname}{Assumption}
\providecommand{\lemmaname}{Lemma}
\providecommand{\remarkname}{Remark}
\providecommand{\theoremname}{Theorem}

\begin{document}
\title{Regression Model Selection Under General Conditions}
\author{Amaze Lusompa \thanks{I thank Sai Avinash, Jason P. Brown, Òscar Jordà, and Francisco Scott
for helpful comments, discussions, and suggestions. I thank Johnson
Oliyide for excellent research assistance. The views expressed are
those of the author and do not necessarily reflect the positions of
the Federal Reserve Bank of Kansas City or the Federal Reserve System.}\\
Federal Reserve Bank of Kansas City\\
}
\maketitle
\begin{abstract}
Model selection criteria are one of the most important tools in statistics.
Proofs showing a model selection criterion is asymptotically optimal
are tailored to the type of model (linear regression, quantile regression,
penalized regression, etc.), the estimation method (linear smoothers,
maximum likelihood, generalized method of moments (GMM), etc.), the
type of data (i.i.d., dependent, high dimensional, etc.), and the
type of model selection criterion. Moreover, assumptions are often
restrictive and unrealistic making it a slow and winding process for
researchers to determine if a model selection criterion is selecting
an optimal model. This paper provides general proofs showing asymptotic
optimality for a wide range of model selection criteria under general
conditions. This paper not only asymptotically justifies model selection
criteria for most situations, but it also unifies and extends a range
of previously disparate results. 
\end{abstract}
\pagebreak{}

\section*{1\quad{}Introduction}

Model selection criteria are one of the most important tools in statistics
and are used for a variety of things such as model selection, model
comparison, model averaging, or to evaluate the forecasting ability
of models. Model selection criteria are used in a wide range of models
from parametric, semiparametric, to nonparametric. Though model selection
criteria have been around for decades and are a lynchpin in statistics,
the behavior of these methods are complex and not fully understood
\citep{Bates2024}. Asymptotic optimality of these criteria (their
ability to select the true model or the best approximating model if
the true model is not in the set), have only been shown for certain
types of models in a limited number of situations. 

Model selection criteria can be broken down into three major classes:
Cross Validation (CV), information criteria, and pseudo-out-of-sample
forecasting.\footnote{Bayesian model selection via the marginal likelihood is a class not
discussed in this paper in part because Bayesian and frequentist notions
of model selection do not necessarily agree \citep{Moreno2015} (see
\cite{Chib2016} and references therein for more on model selection
consistency in the Bayesian case). It is common for forecasting models
estimated using Bayesian methods to use frequentist model selection
criteria such as pseudo-out-of-sample forecasting.}\textsuperscript{,}\footnote{General-to-specific and specific-to-general testing are not a classes
discussed in this paper. Though they have some popularity, unlike
the other major classes they lack wide applicability (e.g. you cannot
use it to find the optimal tuning parameters or things not directly
involving variable selection).} Though CV is considered by many to be the gold standard for model
selection, asymptotic optimality results for CV have been limited.
\cite{Li1987} shows that leave-one-out Cross-Validation (LOO CV)
is asymptotically optimal for a class of linear smoother models under
the assumption that the data are i.i.d. \cite{Andrews1991c} extends
the results to handle regressions with independent but heteroskedastic
residuals. Extending these results to time series, however, has been
spotty at best. Results for CV either assume strict exogeneity, which
is well known to be an unrealistic assumption for time series \citep{Stock2007},
or if they allow for lagged dependent variables, the assumptions for
these models end up being highly restrictive such as requiring the
residuals to be i.i.d. which rules out models with heteroskedasticity
of any kind as well as models with autocorrelated residuals (see for
example \cite{Zhang2013,Sun2021}).\footnote{This rules out most time series models which include but are not limited
to: Local Projections and direct forecast regressions \citep{Jorda2005}
and Vector Autoregressions or autoregressive distributed lag models
with heteroskedasticity of any kind.} As noted in \cite{Hansen2012}, extending the proofs to time-series
regression would be quite challenging.\footnote{To my knowledge there have been scant results for panel data with
\cite{Gao2016,Yu2025} being exceptions.} Moving away from linear smoothers, assumptions are even more restrictive
and things covered are not that broad.\footnote{\cite{Chetverikov2021} show the validity of $k$-fold cross validation
for Lasso in high dimensions assuming the errors are i.i.d. and neither
sub-Gaussian nor sub-exponential. As noted in their paper, the results
do not cover LOO CV.} Popular methods that are not covered at all or only covered under
restrictive conditions include but are not limited to: Lasso and its
variants (adaptive lasso, elastic net, square root lasso, etc.), bridge
estimators more generally, least angle regression (LARS) and other
stepwise methods, quantile regressions, partially linear regression,
nonlinear regressions, generalized linear models, neural networks,
regression trees, ensemble methods, and other machine learning methods. 

Asymptotic optimality of information criteria such as AIC and BIC
have typically only been shown under more stringent assumptions such
as i.i.d. data or i.i.d. residuals (see for example \cite{Shao1997,Claeskens2016,Ding2018}
for extensive reviews). An exception is \cite{Sin1996} who consider
information criteria for dependent processes, but they do not account
for nonparametric or time-varying parameter models. The paper also
assumes restrictive conditions such as the models being finite dimensional
as well as continuous differentiability of the log quasi-likelihood
function which rules out: quantile regressions, many penalized regressions
(such as lasso, its variants, most bridge estimators, and trend filtering),
robust regressions (Huber loss functions), LARS and other stepwise
methods, and regression trees just to name a few. \cite{Sin1996}
acknowledge that generalizing their results to a wider class of data
generating processes and models represents an interesting and challenging
area for future research, but general progress has not been made.\footnote{\cite{Sin1996} show asymptotic optimality by showing the model selection
criteria choose the model with lowest Kullback-Leibler divergence
from the data generating process. In this paper asymptotic optimality
is demonstrated by showing model selection criteria choose the model
with the lowest integrated mean squared error. \cite{Sin1996} also
prove strong consistency of certain information criteria. Strong consistency
arguments are not attempted in this paper. See the discussion in section
4 for the reasoning.} 

Real-time out of sample forecasting exercises (also known as pseudo-out-of-sample
forecasting) in economics and finance are regarded by many researchers
as the \textquotedblleft ultimate test of a forecasting model\textquotedblright{}
in time series (see \cite{Stock2007} page 571).\footnote{This is particularly the case in economics and finance. For example,
central banks routinely evaluate their forecasting models using pseudo-out-of-sample
exercises \citep{Faust2013}, and asset pricing studies typically
require pseudo-out-of-sample evidence for predictability claims to
be taken seriously (see for example \cite{Welch2008}). The preference
for pseudo-out-of-sample testing reflects the high cost of forecast
errors in policy and investment decisions.} These methods are the standard way to choose/evaluate forecast models.
Despite the influx of several new and more complicated methods over
the past few decades (e.g. penalized regression methods, model averaging,
machine learning methods, etc.), asymptotic optimality of pseudo-out-of-sample
forecast methods such as rolling window and recursive methods, are
generally limited to least square estimators or standard time-varying
parameter models (see \cite{Rossi2021} and references therein). 

Currently, researchers spend entire papers or substantial parts of
papers showing asymptotic optimality of a selection procedure, if
it is shown at all, for a limited class of estimators under restrictive
assumptions. In this paper, I derive the asymptotic optimality of
model selection criteria under fairly general conditions. The proofs
do not rely on a specific estimation method and encompass a wide array
of data generating processes. Section 2 reviews CV, information criteria,
and pseudo-out-of-sample forecasting. Section 3 presents the main
proofs for CV, information criteria, and pseudo-out-of-sample forecasting.
Section 4 provides implications and a broader discussion of results
in the literature. Section 5 concludes.

\textbf{Some notation}: $\xrightarrow{p}$ is converges in probability,
$\parallel\cdot\parallel$ is the Frobenius/Euclidean norm, $O_{p}(\cdot)$
is big O probability notation, and $o_{p}(\cdot)$ is little o probability
notation.

\section*{2\quad{}Preliminaries for Cross Validation, Information Criteria,
and Pseudo-Out-Of-Sample Forecasting}

Let the true model be
\[
y_{i}=\mu_{i}+\varepsilon_{i}\;for\;i=1,2,\ldots,T,
\]
where $y_{i}$ is the dependent variable, $\mu_{i}$ is the conditional
mean, and $\varepsilon_{i}$ is the residual. In the special case
of linear regressions
\[
y_{i}=x_{i}\beta+\varepsilon_{i}\;for\;i=1,2,\ldots,T,
\]
where $x_{i}$ is a $1\times p$ vector of regressors, $\beta$ is
a $p\times1$ vector of regression coefficients. The Leave-one-out
(LOO) residual for a regression model is calculated by estimating
$\beta$ for all but one observation, and predicting the residual
for the left out observation. More formally for observation $i$,
\[
\widetilde{\beta}_{-i}=(\sum_{j\neq i}x_{j}'x_{j})^{-1}\sum_{j\neq i}x_{j}'y_{j},\quad\widetilde{\mu}_{-i}=x_{i}\widetilde{\beta}_{-i},\quad\widetilde{\varepsilon}_{i}=y_{i}-\widetilde{\mu}_{-i},
\]
where $\widetilde{\beta}_{-i}$ is the leave $i$ out OLS estimate
of $\beta$, $\widetilde{\mu}_{-i}$ is the leave $i$ out estimate
of $\mu_{i}$, and $\widetilde{\varepsilon}_{i}$ is the leave $i$
out estimate of $\varepsilon_{i}$. LOO CV is calculated by finding
the model that minimizes the mean-squared error (MSE) of the LOO residuals,
that is, by finding the model that minimizes $T^{-1}\sum_{i=1}^{N}\widetilde{\varepsilon}_{i}^{2}$,
though other loss functions such as the absolute value can be used. 

LOO CV was developed and is mainly used for independent data. It is
well known that LOO CV can perform poorly in finite samples when used
for dependent data due to the residuals and regressors being autocorrelated.
This can lead to overfitting or underfitting \citep{Opsomer2001,Arlot2010,Hansen2010}.
One solution is to use block CV methods \citep{Burman1994}. Block
methods work the same way as the LOO method, but instead of leaving
out just observation $i$, one would additionally leave out $h$ observations
on both sides of $i$ to break up the dependence. More formally, one
would leave out observations $i-h,i-h+1,\ldots,i,\ldots i+h-1,i+h$,
so for $h$-block CV via OLS estimation
\[
\widetilde{\beta}_{-i}=(\sum_{j\neq i-h:i+h}x_{j}'x_{j})^{-1}\sum_{j\neq i-h:i+h}x_{j}'y_{j}.
\]
Though fairly popular, the block methods have not been asymptotically
justified.

Another major strand of model selection criteria are information criteria.
They generally take the form of 
\[
log\Bigg(\frac{1}{T}\sum_{i=1}^{T}\{y_{i}-\hat{\mu}_{i}\}^{2}\Bigg)+\frac{\lambda_{T}p}{T}
\]
where in the standard linear regression case $\hat{\mu}_{i}=x_{i}\hat{\beta}$,
$\hat{\beta}$ is the estimate of $\beta$ based on all of the data,
$p$ is the number of predictors and $\lambda_{T}>0$ is the penalty
coefficient. For the AIC $\lambda_{T}=2$, for the BIC $\lambda_{T}=log(T)$,
for HQIC $\lambda_{T}=clog(log(T))$ where $c>2$, for the RIC $\lambda_{T}=2log(p)$,
etc. Models are chosen by choosing the model that minimizes the information
criteria of choice.

The last strand of model selection criteria we will discuss is pseudo-out-of-sample
forecasting. The two most popular methods are rolling window forecast
and recursive forecasts \citep{Clark2009a}.\footnote{In the least squares case, recursive forecasts are also known as predictive
least squares \citep{Wei1992}.} In the standard linear regression case estimated using OLS, rolling
window forecasts estimate a model on the last $R$ observations so
that 
\[
\ddddot{\beta_{i}}=(\sum_{j=i-R}^{i-1}x_{j}'x_{j})^{-1}\sum_{j=i-R}^{i-1}x_{j}'y_{j},\quad\ddddot{\mu_{i}}=x_{i}\ddddot{\beta_{i}},\quad and\quad\ddddot{\varepsilon_{i}}=y_{i}-\ddddot{\mu_{i}}.
\]
Recursive models alternatively estimate the model based on all observations
up to that point with 
\[
\ddddot{\beta_{i}}=(\sum_{j=1}^{i-1}x_{j}'x_{j})^{-1}\sum_{j=1}^{i-1}x_{j}'y_{j}.
\]
The optimal rolling window and recursive forecast models are chosen
by choosing the model that minimizes the MSE of the forecasted residuals.
Note that even though the above examples estimate the standard linear
regression model using OLS, the proofs of asymptotic optimality in
this paper are general and can handle most types of regression models
and estimation methods.

Ideally, the goal of these methods is to minimize the integrated mean
squared error (IMSE). The IMSE is a distance or loss function between
the true conditional mean and the estimated conditional mean for model
$\alpha$ with models and is a measure of how well the conditional
mean for model $\alpha$ approximates the true conditional mean. A
model with a smaller IMSE means that model is closer to the truth
in the squared error loss sense. Define $\mathcal{A}{}_{T}$ as the
set of all models being compared, and $\alpha$ is the model index.\footnote{Note that the models being compared may be individual models or combination/model
averaged models.} The IMSE is defined as
\[
\widetilde{L}_{T}(\alpha)=\frac{1}{T}\sum_{i=1}^{T}(\mu_{i}-\widetilde{\mu}_{-i}(\alpha))^{2}\qquad and\qquad L_{T}(\alpha)=\frac{1}{T}\sum_{i=1}^{T}(\mu_{i}-\hat{\mu}_{i}(\alpha))^{2}
\]
for leave out and full sample estimates respectively. For recursive
and rolling window estimates, the IMSE is define as 
\[
\ddddot{L}_{T}(\alpha)=\frac{1}{T-t_{0}}\sum_{i=t_{0}+1}^{T}(\mu_{i}-\ddddot{\mu_{i}}(\alpha))^{2}\qquad and\qquad\ddddot{L}_{T}(\alpha)=\frac{1}{T-R}\sum_{i=R+1}^{T}(\mu_{i}-\ddddot{\mu_{i}}(\alpha))^{2},
\]
where $t_{0}$ is the minimum number of observations needed so $\ddddot{\mu_{i}}(\alpha)$
is uniquely defined for all models and all observations, $i$.\footnote{For least squares models, $t_{0}$ would be the max $p(\alpha)$ for
all models in the recursive case. In the rolling window case, $R$
would need to be large enough that $\ddddot{\mu_{i}}(\alpha)$ is
uniquely defined.} To show asymptotic optimality, I follow the standard in the literature
and show asymptotically optimality in the sense that 
\[
\frac{L_{T}(\widehat{\alpha})}{\underset{\alpha\in\mathcal{A}_{T}}{inf}L_{T}(\alpha)}\overset{p}{\rightarrow}1
\]
where $\hat{\alpha}$ is the model in the set of $\mathcal{A}{}_{T}$
selected by CV, information criteria, or pseudo out sample forecasting
methods. Intuitively, the above formula says that as the sample size
tends toward infinity, the probability of the model selection procedure
choosing the model with the smallest IMSE converges to 1.\footnote{Note that asymptotic optimality is shown in terms of the full sample
IMSE.}

\section*{3\quad{}Optimality of Model Selection Under General Conditions}

This section shows the proofs for asymptotic optimality. Section 3.1
shows the proofs for CV, section 3.2. for information criteria, and
section 3.3. for pseudo-out-of-sample forecasting.

\subsection*{3.1\quad{}Cross-Validation}

I use the following assumptions to show asymptotic optimality for
CV:
\begin{assumption}
The following conditions are satisfied $\forall\alpha\in\mathcal{A}_{T}$:

(a) For all i, $E(\mu_{i}\varepsilon_{i})=0$, $E(x_{i}(\alpha)\varepsilon_{i})=0$,
$E(\varepsilon_{i}^{2})<\infty$, and each element in $x_{i}(\alpha)\varepsilon_{i}$
has finite second moments.

(b) $\parallel\hat{\mu}_{i}(\alpha)-\widetilde{\mu}_{-i}(\alpha)\parallel\overset{p}{\rightarrow}0$
and $\parallel\hat{\mu}_{i}(\alpha)-\widetilde{\mu}_{-i}(\alpha)\parallel<\infty$$\enspace\forall i$.

(c) $\{\mu_{i}\varepsilon_{i}\}_{i=1}^{T}$ satisfies conditions for
a law of large numbers.

(d) $\{x_{i}(\alpha)\varepsilon_{i}\}_{i=1}^{T}$ satisfies conditions
for a law of large numbers.

(e)The dimension of $x_{i}(\alpha)$, $p(\alpha)$, can grow with
the sample size, but it diverges at a slower rate than the rate of
convergence for the applicable law of large numbers. 

(f)$\underset{\alpha\in\mathcal{A}_{T}}{sup}|\widetilde{L}_{T}(\alpha)/L_{T}(\alpha)-1|\overset{p}{\rightarrow}0$

(g) Each model $\alpha$ can be written as a linear regression $y_{i}=x_{i}(\alpha)\beta(\alpha)+\varepsilon_{i}(\alpha)$
where $\mu_{i}(\alpha)=x_{i}(\alpha)\beta(\alpha)$. 

(h) Either all of the models in the set $\mathcal{A}_{T}$ are misspecified,
or at most one model in the set is true.
\end{assumption}
First note that $\varepsilon_{i}$ is the true residual and $\varepsilon_{i}(\alpha)$
is the residual for model $\alpha$. Assumption 1(a) should be an
uncontroversial assumption since it is made up of standard first order
moment conditions and requires that the second moments for the true
residual and the regression score exist and are bounded. The first
two parts follow directly from the exogeneity assumption, $E(\varepsilon_{i}|x_{i})=0$,
and depending on the models being considered could be strictly (past,
present, and future) exogenous, past and present exogenous, and past
exogenous. Most time series models would be past or past and present
exogenous, while non-time series models are mostly assumed to be strictly
exogenous. Note that the dimension of $x_{i}(\alpha)\varepsilon_{i}$
is $1\times p(\alpha)$. The dimensions are suppressed in the proofs
to make the notation simpler, but the dimensions are taken into account
in the proofs. Assumption 1(a) can even apply to quantile regressions,
though it does not apply to quantile regressions in all situations
(see discussions \cite{Machado2019}).\footnote{In the context of this paper, for quantile regressions one would want
to minimize the IMSE between the conditional quantile for model $\alpha$
and the true conditional quantile. Other loss functions have been
used in the literature for quantile regressions (see for example \cite{Lu2015}). } Note that 1(a) would apply even in the case of omitted variable bias
or classical measurement error since $\varepsilon_{i}$ is the true
residual. If, however, there is simultaneity, 1(a) would not apply
unless one were to first instrument the endogenous variable(s).

The first part of Assumption 1(b) is a more general version of the
assumption made in \cite{Li1987} and \cite{Hansen2012} that the
leverage values for each model $\alpha$ dissipate to zero as $T\rightarrow\infty$.
The second part is just a uniform boundness assumption which rules
out razors edge cases where, for example, the difference grows as
a function of the sample size. As argued in \cite{Li1987}, it should
be the case that the impact of leaving out an observation has on an
estimate should dissipate asymptotically for any reasonable estimator.\footnote{This would clearly apply to M-estimators that satisfy standard regularity
conditions and would apply more generally to extremum estimators that
satisfy certain regularity conditions (see for example Theorem 2.1
in \cite{Newey1994}). It is also consistent with convergence of misspecified
models (see for example \cite{White1981,White1982,White1994} and
references therein).} Note that in the case of $h$-block CV, one can allow $h$ to grow
with the sample size, but to ensure that the condition is satisfied,
one can let $\frac{h}{T}\rightarrow0$ as $h,T\rightarrow\infty$.\footnote{This is standard in block bootstrapping \citep{Lahiri2003}. It may
also be possible to set $h$ as a constant fraction of the sample
size.}

Assumption 1(c) and 1(d) are standard model assumptions that are used
when proving consistency. There are a wide range of law of large numbers
(LLN) for different data including: i.i.d., independent but not identically
distributed, stationary, mixing, mixingale, near-epoch dependent,
and locally stationary just to name a few. For a more exhaustive list
for different types of cross sectional and times series data, see
for example chapter 3 in \cite{White2000} or part 4 in \cite{Davidson2022}.
Assumption 1(e) places a restriction on the the number of variables
that can be included in a model. For example, many law of large numbers
converge at rate $O_{p}(T^{-1/2})$ for each element in the vector,
so assuming $\frac{p(\alpha)^{2}}{T}\rightarrow0$, as $p(\alpha),T\rightarrow\infty$
would be sufficient in many cases. Assuming that $\frac{p(\alpha)^{2}}{T}\rightarrow0$
or $\frac{p(\alpha)^{3}}{T}\rightarrow0$ is standard in infinite
dimensional time series models, infinite dimensional semiparametric
models, as well as high dimensional models \citep{Lutkepohl2005,Chen2007,Fan2004}.\footnote{Note that high-dimensional models with a diverging number of parameters
such that $\frac{p(\alpha)}{T}\in(0,\infty)$ as $p(\alpha),T\rightarrow\infty$
(which are sometimes called proportional asymptotics) are excluded.} This assumption also rules out the razors edge case of perfectly
predicting $y_{i}$ by including as many variables as there are observations.

Assumption 1(f) is standard in the literature, and can be imposed
for example by assuming that either $\parallel\mu_{i}\parallel$ or
$\parallel\mu_{i}-\hat{\mu}_{i}(\alpha)\parallel$ are bounded above
for all $i$ and $\alpha$ in addition to Assumption 1(b) (see Lemma
1 in the appendix). Assumption 1(g) assumes that the candidate models
can be written as linear regressions and is standard assumption in
model selection.\footnote{See for example \cite{White1981,White1982,White1994} and references
therein for more about pseudo-true models.} Though 1(g) is set up for variable/specification selection, the proofs
also apply to bandwidth selection or other tuning parameters where
$\alpha$ denotes the model for a specific tuning parameter. Assumption
1(g) is a broad assumption that includes a broad range of regression
models, however, it excludes nonparametric kernel regression models
(e.g. Nadaraya-Watson or Gasser-Muller types) as well as time-varying
parameter models. These models will require modifications to the proof
and will be addressed in Theorem 2.

Assumption 1(h) is a standard assumption in the literature (see \cite{Li1987,Shao1997,Hansen2012}
and references therein). It is consistent with the statistical adage
from George Box that ``All models are wrong'' or ``All models are
wrong, but some are useful'' as a way to point out that in science,
all models are approximate, so we should not believe the true model
is in the set. This assumption is satisfied for infinite dimensional
models since the true model cannot be in the set. For finite dimensional
models, people do not think the true model is in the set anyway.\footnote{If one includes multiple true models in the set, the model selection
procedures will still select a true model since asymptotically the
model selection procedure will choose a model whose IMSE converges
to zero, though it may not select the most parsimonious version of
the true model in the set whose IMSE will converge to faster than
the other true models \citep{Shao1997}. So it is essentially a razors
edge case where a true model is selected (and the true model will
have a IMSE which converges to zero), but the true model that is selected
does not necessarily have the IMSE which converges to zero the fastest,
so the model selection procedure may not asymptotic optimality in
this case.}
\begin{thm}
Under Assumptions 1 $\hat{\alpha}$ chosen by LOO CV or $h$-block
CV is asymptotically optimal.
\end{thm}
\begin{proof}
Note that $\widetilde{\varepsilon}_{-i}(\alpha)=y_{i}-\widetilde{\mu}_{-i}(\alpha)=\mu_{i}+\varepsilon_{i}-\widetilde{\mu}_{-i}(\alpha)$.
Using simple algebra, it can be shown that the LOO (or the $h$-block
equivalent) squared residual can be decomposed as
\[
\widetilde{\varepsilon}_{-i}^{2}(\alpha)=\varepsilon_{i}^{2}+2(\mu_{i}-\widetilde{\mu}_{-i}(\alpha))\varepsilon_{i}+(\mu_{i}-\widetilde{\mu}_{-i}(\alpha))^{2},
\]
implying that
\[
\frac{1}{T}\sum_{i=1}^{T}\widetilde{\varepsilon}_{-i}^{2}(\alpha)=\frac{1}{T}\sum_{i=1}^{T}\varepsilon_{i}^{2}+\frac{1}{T}\sum_{i=1}^{T}2\mu_{i}\varepsilon_{i}-\frac{1}{T}\sum_{i=1}^{T}2\widetilde{\mu}_{-i}(\alpha)\varepsilon_{i}+\widetilde{L}_{T}(\alpha).
\]
Note that first term, $\frac{1}{T}\sum_{i=1}^{T}\varepsilon_{i}^{2}$,
is the mean squared error of the true residuals and shows up in every
model, so it can be ignored since it will not affect the ordering
of the models. By Assumption 1(f), $\widetilde{L}_{T}(\alpha)$ can
be replaced by $L_{T}(\alpha)$ asymptotically. To show asymptotic
optimality, it is therefore sufficient to show that $\frac{1}{T}\sum_{i=1}^{T}2\mu_{i}\varepsilon_{i}\overset{p}{\rightarrow}0$
and $\frac{1}{T}\sum_{i=1}^{T}2\widetilde{\mu}_{-i}(\alpha)\varepsilon_{i}\overset{p}{\rightarrow}0$
since if this is the case, the LOO MSE of the different models only
would only differ because of their model specific IMSE.\footnote{Technically since $\frac{1}{T}\sum_{i=1}^{T}2\mu_{i}\varepsilon_{i}$
shows up in every model, it can be ignored as well.} Therefore, choosing the model with the smallest LOO MSE would be
choosing the model with the smallest IMSE. 

$\frac{1}{T}\sum_{i=1}^{T}2\mu_{i}\varepsilon_{i}\overset{p}{\rightarrow}0$
follows directly from the exogeneity condition (Assumption 1(a)) and
an appropriate LLN (Assumption 1(c)). To show that $\frac{1}{T}\sum_{i=1}^{T}2\widetilde{\mu}_{-i}(\alpha)\varepsilon_{i}\overset{p}{\rightarrow}0$,
first note that by 1(g) $\widetilde{\mu}_{-i}(\alpha)=x_{i}(\alpha)\widetilde{\beta}_{-i}(\alpha)$.
So 
\[
\frac{1}{T}\sum_{i=1}^{T}2\widetilde{\mu}_{-i}(\alpha)\varepsilon_{i}=\frac{1}{T}\sum_{i=1}^{T}2x_{i}(\alpha)\hat{\beta}(\alpha)\varepsilon_{i}+[\frac{1}{T}\sum_{i=1}^{T}2x_{i}(\alpha)(\widetilde{\beta}_{-i}(\alpha)-\hat{\beta}(\alpha))\varepsilon_{i}]
\]
\[
=\underbrace{[\frac{1}{T}\sum_{i=1}^{T}2x_{i}(\alpha)\varepsilon_{i}]\hat{\beta}(\alpha)}_{\overset{p}{\rightarrow}0}+\underbrace{[\frac{1}{T}\sum_{i=1}^{T}2x_{i}(\alpha)(\widetilde{\beta}_{-i}(\alpha)-\hat{\beta}(\alpha))\varepsilon_{i}]}_{\overset{p}{\rightarrow}0}.
\]
To understand why first term converges to zero, note that $[\frac{1}{T}\sum_{i=1}^{T}2x_{i}(\alpha)\varepsilon_{i}]$
converges in probability to a zero vector by the exogeneity condition
(Assumption 1(a)) and an appropriate LLN for $x_{i}(\alpha)\varepsilon_{i}$
(Assumption 1(d)). In the finite dimension case that is enough for
the first term to converge. If the dimension $p(\alpha)$ grows with
the sample size, by Assumption 1(e) $p(\alpha)$ grows at a slower
rate than the appropriate LLN so the entire first term converges to
zero in that case as well.\footnote{If for example the LLN converges at the standard rate of $O_{p}(T^{-1/2})$,
then $\parallel[\frac{1}{T}\sum_{i=1}^{T}2x_{i}(\alpha)\varepsilon_{i}]\hat{\beta}(\alpha)\parallel=O_{p}(\frac{p(\alpha)}{\sqrt{T}})$
and as long as $p(\alpha)$ grows slower than $\sqrt{T}$, there will
be convergence to zero in probability.} To show why the second term converges to zero, it is sufficient to
show that
\[
E\parallel\frac{1}{T}\sum_{i=1}^{T}2x_{i}(\alpha)(\widetilde{\beta}_{-i}(\alpha)-\hat{\beta}(\alpha))\varepsilon_{i}\parallel\xrightarrow{p}0
\]
since convergence in mean implies convergence in probability. Note
that
\[
E\parallel\frac{1}{T}\sum_{i=1}^{T}2x_{i}(\alpha)(\widetilde{\beta}_{-i}(\alpha)-\hat{\beta}(\alpha))\varepsilon_{i}\parallel\leq\frac{1}{T}\sum_{i=1}^{T}E\parallel2x_{i}(\alpha)(\widetilde{\beta}_{-i}(\alpha)-\hat{\beta}(\alpha))\varepsilon_{i}\parallel
\]
\[
\leq\frac{1}{T}\sum_{i=1}^{T}\underbrace{(E\parallel2x_{i}(\alpha)(\widetilde{\beta}_{-i}(\alpha)-\hat{\beta}(\alpha))\parallel^{2})^{1/2}}_{\overset{p}{\rightarrow}0}\underbrace{(E(\varepsilon_{i}^{2}))^{1/2}}_{bounded}\leq constant\frac{1}{T}\sum_{i=1}^{T}o_{p}(1)=o_{p}(1)
\]
due to the triangle inequality, the Cauchy-Schwarz inequality, Assumption
1(a), and Assumption 1(b).
\end{proof}
\textcompwordmark{}

As was just shown, setting up the problem by decomposing the LOO MSE
in terms of the MSE of the true residuals, the IMSE, and what are
essentially standard first order moment conditions, allows us to leverage
standard statistical results without having to specify the type of
estimation method and easily allows for generality for the type of
model and type of data. The above proof also shows $h$-block CV is
asymptotically optimal, which to my knowledge, no one has actually
shown \citep{Burman1994,Racine2000}. These proofs could also be applied
to $k$-fold CV, but it would require assumptions on the fold size
similar to $h$-block CV (e.g. if you have $k$ folds with $h$ observations
in each fold, $h$ can grow with the sample size but $\frac{h}{T}\rightarrow0$).

As mentioned earlier, Theorem 1 does not allow for observation dependent
coefficients, so it does not account for time-varying parameter models
or Nadaraya-Watson type nonparametric models. To address this shortcoming,
we need a slightly different set of assumptions.
\begin{assumption}
Assumptions 1 holds except for 1(e) and 1(g). In addition, assume
either conditions (a) or (b) holds (conditions a and b are listed
below). Lastly, condition (c) holds (condition c is listed below).

(a) For all $\alpha\in\mathcal{A}_{T}$, $y_{i}=x_{i}(\alpha)\beta(\alpha,x_{i}(\alpha))+\varepsilon_{i}(\alpha)$
where $\mu_{i}(\alpha)=x_{i}(\alpha)\beta(\alpha,x_{i}(\alpha))$
for the nonparametric kernel regression case. $\beta(\alpha,x_{i}(\alpha))$
is continuous in $x_{i}(\alpha)$.

(b) For all $\alpha\in\mathcal{A}_{T}$, $y_{i}=x_{i}(\alpha)\beta_{i}(\alpha)+\varepsilon_{i}(\alpha)$
where $\mu_{i}(\alpha)=x_{i}(\alpha)\beta_{i}(\alpha)$ in the time-varying
parameter case. $\beta_{i}(\alpha)$ is continuous in $i$ in the
infill asymptotic sense.

(c) The data can be divided into $m$ blocks with $\ell$ observations
in each block. These blocks can be constructed in such a way that
$\ell$ goes to infinity, but $\hat{\beta}_{i}(\alpha)=\hat{\beta}_{r}(\alpha)+o_{p}(1)$
or $\hat{\beta}_{}(\alpha,x_{i})=\hat{\beta}_{}(\alpha,x_{r})+o_{p}(1)$
for $i\neq r$ within the same block.
\end{assumption}
Assumption 2 is essentially the same as Assumption 1, but there are
two differences. The first difference is it allows the model parameters
to vary with the observation or covariate. For time-varying parameter
models, $\beta_{i}(\alpha)$ is a continuous function in time and
follows standard infill asymptotics arguments \citep{Robinson1989,Cai2007,Chen2012,Dahlhaus2012}.
Note that for time-varying parameter case, the continuous function
can be deterministic or stochastic.\footnote{The standard in nonparametric time-varying parameter estimation is
to use nonparametric kernels \citep{Robinson1989,Cai2007,Dahlhaus2012},
and it is typically assumed that the function is deterministic and
differentiable. It has been shown that kernel methods can be used
for certain stochastic processes \citep{Giraitis2014,Giraitis2021}.
As long as the continuous functions can be written as an infinite
order of orthogonal basis functions (whether it be by the Stone-Weierstrass
Theorem or Karhunen-Loève Theorem), the time-varying parameters can
alternatively be estimated using basis function approximations \citep{Huang2002,Huang2004}.} The continuity assumption for nonparametric kernel regressions is
also standard \citep{Fan1996}. Assumption 2(c) just says as the sample
size grows, the difference between estimates in the same block should
decrease asymptotically. Note that the construction of the $m$ blocks
of size $\ell$ is in line with consistency assumptions and intuition
behind nonparametric kernel regressions and infill asymptotics of
time-varying parameter models. In the case where the estimators are
consistent/pseudo consistent, this condition is satisfied.\footnote{It is fairly standard in nonparametric regression literature to assume
that the class of estimators being compared are consistent/pseudo
consistent (see for example \cite{Li1985,Hardle1994,Hart1994,Leung2005,Sun2021}
and references therein).} There are also cases where this condition can be satisfied even if
the estimators are not consistent/pseudo consistent (e.g. rolling
window regression where the window is a constant fraction of the sample
size or more generally a kernel regression where the bandwidth does
not go to zero fast enough to satisfy the consistency condition).

The second difference from Assumption 1 is that the dimension, $p(\alpha)$,
is not growing with the sample size. This is standard in the time-varying
parameter and nonparametric kernel regression literature. As will
be seen in the proofs, the proofs may be able to handle $p(\alpha)$
growing with the sample size, but it require being specific about
the rates at which $p(\alpha)$ and $\ell$ grow.
\begin{thm}
Under Assumptions 2, $\hat{\alpha}$ chosen by LOO CV or $h$-block
CV is asymptotically optimal.
\end{thm}
\begin{proof}
The proof is the same as Theorem 1, except for the proof that $\frac{1}{T}\sum_{i=1}^{T}2\widetilde{\mu}_{-i}(\alpha)\varepsilon_{i}$
converges in probability to 0. The proof will be written for time-varying
parameter model case, but it also applies to the nonparametric kernel
regression. It is written in terms of the time-varying parameter model
for convenience since observations are naturally ordered in regards
to time. To show $\frac{1}{T}\sum_{i=1}^{T}2\widetilde{\mu}_{-i}(\alpha)\varepsilon_{i}$
converges in probability to 0, first divide the data into $m$ blocks
with $\ell$ observations in each block. The blocks are constructed
in such a way that $\ell$ goes to infinity, but $\ell$ is small
enough such that for any two time-varying parameter vectors in the
same block $\hat{\beta}_{i}(\alpha)=\hat{\beta}_{r}(\alpha)+o_{p}(1)$.
For any one block we have 
\[
\frac{1}{\ell}\sum_{i=1}^{\ell}2\widetilde{\mu}_{-i}(\alpha)\varepsilon_{i}=\frac{1}{\ell}\sum_{i=1}^{\ell}2x_{i}(\alpha)\hat{\beta}_{i}(\alpha)\varepsilon_{i}+[\frac{1}{\ell}\sum_{i=1}^{\ell}2x_{i}(\alpha)(\widetilde{\beta}_{-i}(\alpha)-\hat{\beta}_{i}(\alpha))\varepsilon_{i}]
\]
\[
=\underbrace{[\frac{1}{\ell}\sum_{i=1}^{\ell}2x_{i}(\alpha)\varepsilon_{i}]\hat{\beta}_{r}(\alpha)}_{\overset{p}{\rightarrow}0}+\frac{1}{\ell}\sum_{i=1}^{\ell}2x_{i}(\alpha)\varepsilon_{i}o_{p}(1)+\underbrace{[\frac{1}{\ell}\sum_{i=1}^{\ell}2x_{i}(\alpha)(\widetilde{\beta}_{-i}(\alpha)-\hat{\beta}_{i}(\alpha))\varepsilon_{i}]}_{\overset{p}{\rightarrow}0}
\]
where convergence in probability to zero for the first term follows
from the exogeneity condition (Assumption 1(a)) and an appropriate
LLN (Assumption 1(d)). Convergence of the third terms follows the
same argument used in the proof of Theorem 1 and is omitted for brevity.
To show convergence of the second term, note that it is sufficient
to show that 
\[
E\parallel\frac{1}{\ell}\sum_{i=1}^{\ell}2x_{i}(\alpha)\varepsilon_{i}o_{p}(1)\parallel\xrightarrow{p}0.
\]
Similar to an argument used in the proof used in Theorem 1, note that
\[
E\parallel\frac{1}{\ell}\sum_{i=1}^{\ell}2x_{i}(\alpha)\varepsilon_{i}o_{p}(1)\parallel\leq\frac{1}{\ell}\sum_{i=1}^{\ell}\underbrace{(E\parallel2x_{i}(\alpha)\varepsilon_{i}\parallel^{2})^{1/2}}_{O_{p}(p(\alpha))=O_{p}(1)}\underbrace{(E\parallel o_{p}(1)\parallel^{2})^{1/2}}_{\overset{p}{\rightarrow}0}
\]
\[
\leq constant\times\frac{1}{\ell}\sum_{i=1}^{\ell}o_{p}(1)=o_{p}(1)
\]
due to the triangle inequality, the Cauchy-Schwarz inequality, and
Assumption 1(a). Since for each block we have $\parallel\frac{1}{\ell}\sum_{i=1}^{\ell}2\widetilde{\mu}_{-i}(\alpha)\varepsilon_{i}\parallel$
converging to zero in probability, the average over those blocks will
also converge to zero, thus completing the proof.
\end{proof}
\begin{rem}
For the time-varying parameter case, it is assumed that the parameters
map to a continuous function. This may be unrealistic if one thinks
that the time-varying parameter process is not continuous and instead
exhibits structural breaks. The mapping to a continuous function is
not a necessary assumption. In the case of structural breaks, the
proofs still go through as is, but the blocks need to correspond to
break fractions \citep{Andrews1993,Bai1998}. For the non-parametric
kernel regression case, these proofs could also be applied to $k$-fold
CV.\footnote{Again, though 2(a) and 2(b) are set up for variable/specification
selection, these proofs also apply to bandwidth selection or other
tuning parameters where $\alpha$ denotes the model for a specific
tuning parameter.}
\end{rem}

\subsection*{3.2\quad{}Information Criteria}

I use the following assumptions to show asymptotic optimality for
information criteria:
\begin{assumption}
Either Assumption 1 or Assumption 2 is satisfied (except for 1(b)
and 1(f)). In addition, assume the penalty term $\frac{\lambda_{T}p(\alpha)}{T}\rightarrow0$
for all $\alpha\in\mathcal{A}_{T}$
\end{assumption}
Assumptions 1(b) and 1(f) are dropped since information criteria use
all of the data when estimating the residuals unlike CV. The penalty
term tending toward zero asymptotically occurs for most if not all
of the popular information criteria including but not limited to:
AIC, BIC, HQIC, RIC, and Generalized information criteria.
\begin{thm}
Under Assumption 3, $\hat{\alpha}$ chosen by an information criteria
is asymptotically optimal.
\end{thm}
\begin{proof}
The proof follows along similar lines as the proof of Theorems 1 and
2 except there is no need to show $\parallel\frac{1}{T}\sum_{i=1}^{T}2x_{i}(\alpha)(\widetilde{\beta}_{-i}(\alpha)-\hat{\beta}(\alpha))\varepsilon_{i}]\parallel$
or $\parallel[\frac{1}{\ell}\sum_{i=1}^{\ell}2x_{i}(\alpha)(\widetilde{\beta}_{-i}(\alpha)-\hat{\beta}_{i}(\alpha))\varepsilon_{i}]\parallel$
converge to 0 in probability for the constant parameter and time-varying
parameter model/nonparametric model cases, respectively. It is sufficient
to just demonstrate the proof under the constant parameter case and
the time-varying/nonparametric regression cases are omitted for brevity
since they are simpler versions of the proof of Theorem 2. In the
constant parameter case with information criteria, we have
\[
\frac{1}{T}\sum_{i=1}^{T}\hat{\varepsilon}_{i}^{2}(\alpha)=\frac{1}{T}\sum_{i=1}^{T}\varepsilon_{i}^{2}+\frac{1}{T}\sum_{i=1}^{T}2\mu_{i}\varepsilon_{i}-\frac{1}{T}\sum_{i=1}^{T}2\hat{\mu}_{i}(\alpha)\varepsilon_{i}+L_{T}(\alpha).
\]
The penalty has no impact on the ordering of the models asymptotically
since the penalty converges to 0 asymptotically, and the natural log
does not change the ordering of the information criteria for the models
since it is a strictly monotonic function. Therefore to show asymptotic
optimality, it is sufficient to show that $\frac{1}{T}\sum_{i=1}^{T}2\mu_{i}\varepsilon_{i}\overset{p}{\rightarrow}0$
and $\frac{1}{T}\sum_{i=1}^{T}2\hat{\mu}_{i}(\alpha)\varepsilon_{i}\overset{p}{\rightarrow}0$
since by the continuous mapping theorem the only difference between
the information criteria asymptotically would be the model specific
IMSE. The proof $\frac{1}{T}\sum_{i=1}^{T}2\mu_{i}\varepsilon_{i}\overset{p}{\rightarrow}0$
follows same argument used in proof of Theorem 1 and after noting
that $\hat{\mu}_{i}(\alpha)=x_{i}(\alpha)\hat{\beta}(\alpha)$, $\frac{1}{T}\sum_{i=1}^{T}2\hat{\mu}_{i}(\alpha)\varepsilon_{i}\overset{p}{\rightarrow}0$
follows from the argument that $[\frac{1}{T}\sum_{i=1}^{T}2x_{i}(\alpha)\varepsilon_{i}]\hat{\beta}(\alpha)\overset{p}{\rightarrow}0$
shown in Theorem 1.
\end{proof}
So despite AIC, BIC, and other information criteria being derived
under fairly restrictive conditions (i.i.d. data or independent residuals),
they are asymptotically optimal under more general conditions. This
result holds for any penalty that converges to 0 asymptotically, thus
indicating the arbitrariness of the penalty from an asymptotic optimality
standpoint. In finite samples, the chosen model may vary drastically
due to the penalty, but currently it is not clear what the penalty
should be, and it may be the case that there is no best penalty for
all situations.

\subsection*{3.3\quad{}Pseudo-Out-Of-Sample Forecast Methods}

\subsubsection*{3.3.1\quad{}Constant Parameter Case}

To show asymptotic optimality of rolling window and recursive forecasts
in the constant parameter case, we need updated assumptions.
\begin{assumption}
Assume as $R$ and $t_{0}$ grow, Assumption 1 holds. In Assumption
1, $\widetilde{\beta}_{-i}(\alpha)$, $\widetilde{\mu}_{-i}(\alpha)$,
and $\widetilde{L}_{T}(\alpha)$ are replaced by $\ddddot{\beta}_{i}(\alpha)$,
$\ddddot{\mu}_{i}(\alpha)$, and $\ddddot{L}_{T}(\alpha)$, respectively.
In addition, assume

(a) $R$ and $t_{0}$ grow with the sample size but at a slower rate
such that $R,t_{0},T\rightarrow\infty$ but $T-R,T-t_{0}\rightarrow\infty$.
\end{assumption}
Assumption 4 is the rolling window and recursive version of Assumption
1. Assumption 4(a), which requires that the window size, $R$, grows
with the sample size but at a slower rate, is standard in the literature
(see \cite{Inoue2017} and references therein). Letting the initial
starting point for the recursive forecast, $t_{0}$, tend toward infinity
is also common (e.g. \cite{West1996,West1998,Clark2009a}), though
many papers also assume that it is either constant or a constant fraction
of the sample size.\footnote{The proofs for the recursive case can be written under the alternative
assumptions for $t_{0}$, though it would require different proofs,
and it is convenient to assume that it grows with the sample size
as it allows the same proof to be used for both the rolling window
and recursive forecast cases.}
\begin{thm}
Under Assumption 4, $\hat{\alpha}$ chosen by rolling window and recursive
forecast is asymptotically optimal. 
\end{thm}
\begin{proof}
See appendix.
\end{proof}

\subsubsection*{3.3.1\quad{}Time-varying Parameter Case}

To show asymptotic optimality for rolling window and recursive forecasts
in the time-varying parameter case, we use the following assumption:
\begin{assumption}
Assume $R$ and $t_{0}$ grow in such a way that Assumption 2 holds.
In Assumption 2, $\widetilde{\beta}_{-i}(\alpha)$, $\widetilde{\mu}_{-i}(\alpha)$,
and $\widetilde{L}_{T}(\alpha)$ are replaced by $\ddddot{\beta}_{i}(\alpha)$,
$\ddddot{\mu}_{i}(\alpha)$, and $\ddddot{L}_{T}(\alpha)$, respectively.
In addition, assume

(a) $R,t_{0},T\rightarrow\infty$, but $T-R,T-t_{0}\rightarrow\infty$.
\end{assumption}
As long as $R$ and $t_{0}$ grow at suitable rates, the modified
Assumption 2 should be satisfied. In the case of nonparametric kernel
estimation with time-varying parameters, rolling window estimates
are simply one sided kernel estimates of the time-varying parameters
\citep{Inoue2017,Cai2023,Farmer2023}, as opposed to the full sample
estimates which are generally based on two sided kernels. Taking that
into account along with the infill asymptotic framework, the recursive/rolling
window version of Assumption 2(c) is not unrealistic.\footnote{This also makes sense if one uses the basis function approximation
approach to estimating time-varying parameter models \citep{Huang2002,Huang2004}.} Note that the rate at which the block size, $\ell$, grows may be
dependent on the rate at which $R$ and $t_{0}$ grow. Assumption
5(a) is the same as 4(a) and follows the same reasoning.
\begin{thm}
Under Assumption 5, $\hat{\alpha}$ chosen by rolling window and recursive
forecasts is asymptotically optimal. 
\end{thm}
\begin{proof}
See appendix.
\end{proof}

\section*{4\quad{}Implications and Discussion}

Under fairly general conditions, CV, the most popular information
criteria, and the most popular pseudo-out-of-sample forecasting methods
will asymptotically all choose the model with the same conditional
mean, so from an asymptotic optimality standpoint, there is no benefit
from using one procedure versus another. In regards to what model
selection procedure one should use in finite samples, the arguments
in this paper cannot speak to them. I do think the results help highlight
why recommendations for certain model selection procedures are not
as strong as they first appear. 

Many people over the last few decades have advocated the use of model
selection criteria such as BIC and HQIC on the grounds that they are
consistent (see for example \cite{Shao1997,Claeskens2016,Ding2018}
and references therein). What does not get talked about enough is
consistency is a distinction without a difference. Intuitively, consistency
of a model selection procedure simply means that if the most parsimonious
version of the true model is in a set, then the model selection procedure
will choose the most parsimonious version of the true model asymptotically.
What does not get brought up enough is that it could be the case that
in finite samples, the researcher would prefer a criteria that is
not consistent. Consistent model selection criteria tend to choose
a more parsimonious version of the model, which may not be ideal.
There is a host of Monte Carlo evidence indicating this is the case.
For example \cite{Lutkepohl2005} (page 155) shows in a Monte Carlo
where the true model is a vector autoregression (VAR) with two lags
that AIC selected the true model more often than BIC and HQIC, despite
the latter two criteria being consistent. \cite{Burnham2004} shows
in a Monte Carlo that when the true model is denser (has more variables),
AIC does a better job of selecting a model closer to the truth than
BIC, but the results depend on the size of the parameters. \cite{Ng2005}
and \cite{Ng2013} show that neither AIC or BIC dominate, but again,
the results depend on the size of the parameters.

If the true model is not in the set (which would be the case for infinite
dimensional models), a model selection procedure cannot be consistent
since the true model can never be in the set of models being compared.
Even if the true model is finite dimensional, if the true model is
not in the set, then we currently do not have a reason (from the asymptotic
point of view) to prefer one model selection procedure over another.
This is in line with the George Box argument that ``All models are
wrong'' or ``All models are wrong, but some are useful'' as a way
to point out that in science, all models are approximate, so we should
not believe the true model is in the set.\footnote{Many researchers have realized this, and it is one of the reasons
why model averaging has become popular over the past 30 years \citep{Steel2020}. }

There has also been focus on model selection procedures, such as AIC,
being minimax.\footnote{Conditions for minimax estimators are generally derived under i.i.d.
data, and to my knowledge these results have not been extended.}\textsuperscript{,}\footnote{\cite{Yang2005} shows an estimator cannot be strongly consistent
and minimax at the same time.} A model selection procedure being minimax simply means that in the
worst case scenario, the minimax model selection procedure has a smaller
IMSE than a non-minimax procedure, or alternatively, minimax procedures
minimize downside risk. But for any particular situation a researcher
is interested in, a procedure being minimax does not say anything
about which criteria is preferred in finite samples because the situation
may be far from the worst case scenario. Minimax procedures tend to
choose larger models, and consistent procedures tend to choose smaller
models, so depending how accurately and efficiently estimated parameters
are and whether you are more worried about overfitting or underfitting,
you may prefer one procedure over another. In a finite sample, there
is a bias variance tradeoff because while parameters may be biased
when the selected model is too small, the parameter estimates may
be highly inefficient if the estimated model includes too many parameters
\citep{Ng2013,Rossi2021}. Furthermore, as long as $\lambda_{T}$
is a finite constant greater than 1, the information criteria is minimax
\citep{Shao1997,Yang2005}. So even within the class of minimax procedures,
it is not clear which one should be chosen since as shown in Theorem
3, they are all asymptotically optimal.\footnote{The same holds for consistent model selection procedures. A sufficient
condition for consistent model selection procedures require that $\lambda_{T}\rightarrow\infty$
but $\frac{\lambda_{T}p}{T}\rightarrow0$. Since there are an infinite
number of ways to construct a consistent model selection procedures
and they are asymptotically optimal under the conditions stated in
this paper, the literature currently cannot distinguish between procedures
within this class from an asymptotic standpoint.} 

Based off the results in this paper and in the literature, I believe
it is the case that the literature should focus more on the finite
sample properties or maybe use local asymptotics or set up problems
in such a way that the differences in model selection criteria show
up in the limit.\footnote{An important topic of research in the literature that is beyond the
scope of this paper is taking into account the impact model selection
has on subsequent inference (see for example \cite{Leeb2005}).} 

\section*{5\quad{}Conclusion}

This paper provides general proofs showing optimality for a wide range
of model selection criteria under fairly general conditions. This
paper not only asymptotically justifies model selection criteria for
most situations, but it also unifies and extends a range of previously
disparate results. 

The results from this paper should allow researchers to move on from
showing the asymptotic optimality of model selection procedures for
most situations and to potentially focus on showing the theoretical
finite sample properties of these methods. Since it is the case that
the most popular methods end up being asymptotically optimal under
general conditions, the choice of which model selection procedure
to choose from, like most things in statistics, appears to be a finite
sample choice and not an asymptotic one.

\pagebreak\bibliographystyle{chicago}
\bibliography{References}

\pagebreak{}

\subsubsection*{Proof of Theorem 4}
\begin{proof}
The proofs are written for the rolling window case. Note that for
the recursive case, one would just replace $R$ with $t_{0}$. The
proof follows along similar lines as the proof of Theorem 1. Note
that 
\[
\ddddot{\varepsilon}_{i}^{2}(\alpha)=\varepsilon_{i}^{2}+2(\mu_{i}-\ddddot{\mu_{i}}(\alpha))\varepsilon_{i}+(\mu_{i}-\ddddot{\mu_{i}}(\alpha))^{2},
\]
which implies
\[
\frac{1}{T-R}\sum_{i=R+1}^{T}\ddddot{\varepsilon}_{i}^{2}(\alpha)=\frac{1}{T-R}\sum_{i=R+1}^{T}\varepsilon_{i}^{2}+\frac{1}{T-R}\sum_{i=R+1}^{T}2\mu_{i}\varepsilon_{i}-\frac{1}{T-R}\sum_{i=R+1}^{T}2\ddddot{\mu_{i}}(\alpha)\varepsilon_{i}+\ddddot{L}_{T}(\alpha).
\]
To prove asymptotic optimality, it is sufficient to show that $\frac{1}{T-R}\sum_{i=R+1}^{T}2\mu_{i}(\alpha)\varepsilon_{i}\overset{p}{\rightarrow}0$
and $\frac{1}{T-R}\sum_{i=R+1}^{T}2\ddddot{\mu_{i}}(\alpha)\varepsilon_{i}\overset{p}{\rightarrow}0$.
Again, the proof $\frac{1}{T-R}\sum_{i=R+1}^{T}2\mu_{i}(\alpha)\varepsilon_{i}\overset{p}{\rightarrow}0$
follows same arguments used in Theorem 1. To show $\frac{1}{T-R}\sum_{i=R+1}^{T}2\ddddot{\mu_{i}}(\alpha)\varepsilon_{i}\overset{p}{\rightarrow}0$,
note that
\[
\frac{1}{T-R}\sum_{i=R+1}^{T}2\ddddot{\mu_{i}}(\alpha)\varepsilon_{i}=[\frac{1}{T-R}\sum_{i=R+1}^{T}2x_{i}(\alpha)\widehat{\beta}(\alpha)\varepsilon_{i}]+\frac{1}{T-R}\sum_{i=R+1}^{T}[2x_{i}(\alpha)(\ddddot{\beta}_{i}(\alpha)-\widehat{\beta}(\alpha))\varepsilon_{i}]
\]
\[
=\underbrace{[\frac{1}{T-R}\sum_{i=R+1}^{T}2x_{i}(\alpha)\varepsilon_{i}]\widehat{\beta}(\alpha)}_{\overset{p}{\rightarrow}0}+\underbrace{[\frac{1}{T-R}\sum_{i=R+1}^{T}2x_{i}(\alpha)(\ddddot{\beta}_{i}(\alpha)-\widehat{\beta}(\alpha))\varepsilon_{i}]}_{\overset{p}{\rightarrow}0}
\]
where the first term converges to zero by the exogeneity condition
(Assumption 1(a)), an appropriate LLN (Assumption 1(d)), and in the
infinite dimensional case Assumption 1(e). To show that the second
term converges to 0, note that
\[
E\parallel\frac{1}{T-R}\sum_{i=R+1}^{T}2x_{i}(\alpha)(\ddddot{\beta}_{i}(\alpha)-\widehat{\beta}(\alpha))\varepsilon_{i}\parallel\leq\frac{1}{T-R}\sum_{i=R+1}^{T}(E\parallel2x_{i}(\alpha)(\ddddot{\beta}_{i}(\alpha)-\widehat{\beta}(\alpha))\parallel^{2})^{1/2}\underbrace{(E(\varepsilon_{i}^{2}))^{1/2}}_{bounded}
\]
\[
\leq constant\frac{1}{T-R}\sum_{i=R+1}^{T}\underbrace{(E\parallel2x_{i}(\alpha)(\ddddot{\beta}_{i}(\alpha)-\widehat{\beta}(\alpha))\parallel^{2})^{1/2}}_{\overset{p}{\rightarrow}0}=o_{p}(1)
\]
due to the triangle inequality, Cauchy-Schwarz inequality, Assumption
1(a), and Assumption 4(a), thus proving that $\frac{1}{T-t_{0}}\sum_{i=t_{0}+1}^{T}2\ddddot{\mu_{i}}(\alpha)\varepsilon_{i}\overset{p}{\rightarrow}0$. 
\end{proof}

\subsubsection*{Proof of Theorem 5}
\begin{proof}
The proofs are written for the rolling window case. Note that for
the recursive case, one would just replace $R$ with $t_{0}$. Note
that 
\[
\frac{1}{T-R}\sum_{i=R+1}^{T}\ddddot{\varepsilon}_{i}^{2}(\alpha)=\frac{1}{T-R}\sum_{i=R+1}^{T}\varepsilon_{i}^{2}+\frac{1}{T-R}\sum_{i=R+1}^{T}2\mu_{i}\varepsilon_{i}-\frac{1}{T-R}\sum_{i=R+1}^{T}2\ddddot{\mu_{i}}(\alpha)\varepsilon_{i}+\ddddot{L}_{T}(\alpha).
\]
To prove asymptotic optimality, it is sufficient to show that $\frac{1}{T-R}\sum_{i=R+1}^{T}2\mu_{i}\varepsilon_{i}\overset{p}{\rightarrow}0$
and $\frac{1}{T-R}\sum_{i=R+1}^{T}2\ddddot{\mu_{i}}(\alpha)\varepsilon_{i}\overset{p}{\rightarrow}0$.
Again, the proof $\frac{1}{T-R}\sum_{i=R+1}^{T}2\mu_{i}\varepsilon_{i}\overset{p}{\rightarrow}0$
follows same arguments used in Theorem 1. The $\frac{1}{T-R}\sum_{i=R+1}^{T}2\ddddot{\mu_{i}}(\alpha)\varepsilon_{i}\overset{p}{\rightarrow}0$
follows the same argument used in the proof of Theorem 2, except $(\widetilde{\beta}_{-i}(\alpha)-\hat{\beta}_{i}(\alpha))$
is replaced with $(\ddddot{\beta}_{i}(\alpha)_{}-\hat{\beta}_{i}(\alpha))$. 
\end{proof}
\begin{lem}
Assume 1(b) holds. In addition, assume either $\parallel\mu_{i}\parallel$
or $\parallel\mu_{i}-\hat{\mu}_{i}(\alpha)\parallel$ are bounded
above for all $i$ and $\alpha$. Then $|\widetilde{L}_{T}(\alpha)/L_{T}(\alpha)-1|\overset{p}{\rightarrow}0$
for all $\alpha\in A_{T}$.
\end{lem}
\begin{proof}
The proof is shown in two cases. The first under the case $\parallel\mu_{i}\parallel<\infty$
and the second under $\parallel\mu_{i}-\hat{\mu}_{i}(\alpha)\parallel<\infty$.
Under case 1, note that 
\[
\widetilde{L}_{T}(\alpha)-L_{T}(\alpha)=\frac{1}{T}\sum_{i=1}^{T}(\mu_{i}-\widetilde{\mu}_{-i}(\alpha))^{2}-\frac{1}{T}\sum_{i=1}^{T}(\mu_{i}-\hat{\mu}_{i}(\alpha))^{2}
\]
\[
=\frac{1}{T}\sum_{i=1}^{T}2(\mu_{i}\hat{\mu}_{i}(\alpha)-\mu_{i}\widetilde{\mu}_{-i}(\alpha))+\frac{1}{T}\sum_{i=1}^{T}(\widetilde{\mu}_{-i}(\alpha)^{2}-\hat{\mu}_{i}(\alpha)^{2})
\]
\[
=\frac{1}{T}\sum_{i=1}^{T}2(\mu_{i}(\hat{\mu}(\alpha)-\widetilde{\mu}_{-i}(\alpha))+\frac{1}{T}\sum_{i=1}^{T}(\widetilde{\mu}_{-i}(\alpha)^{2}-\hat{\mu}_{i}(\alpha)^{2})
\]
To complete the proof for case 1, it is sufficient to show 
\[
E\parallel\frac{1}{T}\sum_{i=1}^{T}2(\mu_{i}(\hat{\mu}(\alpha)-\widetilde{\mu}_{-i}(\alpha))+\frac{1}{T}\sum_{i=1}^{T}(\widetilde{\mu}_{-i}(\alpha)^{2}-\hat{\mu}_{i}(\alpha)^{2})\parallel\overset{p}{\rightarrow}0.
\]
Note that 
\[
E\parallel\frac{1}{T}\sum_{i=1}^{T}2(\mu_{i}(\hat{\mu}(\alpha)-\widetilde{\mu}_{-i}(\alpha))+\frac{1}{T}\sum_{i=1}^{T}(\widetilde{\mu}_{-i}(\alpha)^{2}-\hat{\mu}_{i}(\alpha)^{2})\parallel
\]
\[
\leq\frac{1}{T}\sum_{i=1}^{T}\underbrace{(E\parallel2\mu_{i}\parallel^{2})^{1/2}}_{bounded}\underbrace{(E\parallel\hat{\mu}(\alpha)-\widetilde{\mu}_{-i}(\alpha)\parallel^{2})^{1/2}}_{\overset{p}{\rightarrow}0}+\frac{1}{T}\sum_{i=1}^{T}\underbrace{E\parallel(\widetilde{\mu}_{-i}(\alpha)^{2}-\hat{\mu}_{i}(\alpha)^{2})\parallel}_{\overset{p}{\rightarrow}0}
\]
\[
=o_{p}(1)
\]
by the triangle inequality, Assumption 1(b), and the assumption of
$\parallel\mu_{i}\parallel$ being bounded for all $i$.

For case 2, note that
\[
(\mu_{i}-\widetilde{\mu}_{-i}(\alpha))^{2}=([\mu_{i}-\hat{\mu}_{i}(\alpha)]+[\hat{\mu}_{i}(\alpha)-\widetilde{\mu}_{-i}(\alpha)])^{2}
\]
\[
=[\mu_{i}-\hat{\mu}_{i}(\alpha)]^{2}+2[\mu_{i}-\hat{\mu}_{i}(\alpha)][\hat{\mu}_{i}(\alpha)-\widetilde{\mu}_{-i}(\alpha)]+[\hat{\mu}_{i}(\alpha)-\widetilde{\mu}_{-i}(\alpha)]{}^{2}.
\]
It follows that 
\[
\widetilde{L}_{T}(\alpha)=L_{T}(\alpha)+\frac{1}{T}\sum_{i=1}^{T}2[\mu_{i}-\hat{\mu}_{i}(\alpha)][\hat{\mu}_{i}(\alpha)-\widetilde{\mu}_{-i}(\alpha)]+\frac{1}{T}\sum_{i=1}^{T}[\hat{\mu}_{i}(\alpha)-\widetilde{\mu}_{-i}(\alpha)]{}^{2}.
\]
To show $|\widetilde{L}_{T}(\alpha)/L_{T}(\alpha)-1|\overset{p}{\rightarrow}0$,
it is sufficient to show that the second and third terms on the right
hand side converge to zero in probability. Convergence of the third
term follows from Assumption 1(b). To show the second term converges
to zero in probability, note that it is sufficient to show that 
\[
E\parallel\frac{1}{T}\sum_{i=1}^{T}2[\mu_{i}-\hat{\mu}_{i}(\alpha)][\hat{\mu}_{i}(\alpha)-\widetilde{\mu}_{-i}(\alpha)]\parallel\xrightarrow{p}0.
\]
Note that 
\[
E\parallel\frac{1}{T}\sum_{i=1}^{T}2[\mu_{i}-\hat{\mu}_{i}(\alpha)][\hat{\mu}_{i}(\alpha)-\widetilde{\mu}_{-i}(\alpha)]\parallel\leq\frac{1}{T}\sum_{i=1}^{T}E\Bigg(\underbrace{\parallel2[\mu_{i}-\hat{\mu}_{i}(\alpha)]\parallel}_{bounded}\parallel[\hat{\mu}_{i}(\alpha)-\widetilde{\mu}_{-i}(\alpha)]\parallel\bigg)
\]
\[
\leq constant\frac{1}{T}\sum_{i=1}^{T}\underbrace{E\Bigg(\parallel[\hat{\mu}_{i}(\alpha)-\widetilde{\mu}_{-i}(\alpha)]\parallel\bigg)}_{\overset{p}{\rightarrow}0}=o_{p}(1),
\]
by the triangle inequality, by the assumption that $\parallel\mu_{i}-\hat{\mu}_{i}(\alpha)\parallel$
is bounded above for all $i$ and $\alpha$, and Assumptions 1(b).
\end{proof}
\begin{rem}
The proof of Lemma 1 also applies for the rolling window and recursive
cases for both the constant parameter and time-varying cases. One
just needs to replace $\widetilde{\beta}_{-i}(\alpha)$, $\widetilde{\mu}_{-i}(\alpha)$,
and $\widetilde{L}_{T}(\alpha)$ with $\ddddot{\beta}_{i}(\alpha)$,
$\ddddot{\mu}_{i}(\alpha)$, and $\ddddot{L}_{T}(\alpha)$, respectively.
One also needs to use the corresponding assumptions in Assumptions
4 and 5.
\end{rem}

\end{document}